\documentclass[11pt, reqno]{amsart}

\usepackage[a4paper,
            top=30mm,bottom=30mm,
            left=36mm,right=36mm]{geometry}

\usepackage{amsmath, amsthm, amsfonts, amssymb, amscd, mathtools, bm} 
\usepackage{xcolor} 
\usepackage[unicode=true]{hyperref} 
\usepackage[shortlabels]{enumitem} 
\usepackage[T1]{fontenc}
\usepackage{empheq}
\usepackage{placeins}

\theoremstyle{plain}
\newtheorem{theorem}{Theorem}[section]

\theoremstyle{definition}

\newtheorem{example}[theorem]{Example}

\theoremstyle{remark}

\numberwithin{equation}{section}

\newcommand{\eps}{\varepsilon}
\newcommand{\teps}{\tilde{\varepsilon}}
\newcommand{\dx}{\mathrm{d} x}
\newcommand{\dz}{\mathrm{d} z}
\newcommand{\rd}{\mathrm{d}}
\newcommand{\N}{\mathbb{N}}
\newcommand{\R}{\mathbb{R}}

\newcommand{\F}{\mathcal{F}}

\title[]{An integral characterization \\ of almost equicontinuity}

\author[]{Nuno J. Alves}
\author[]{Hikmatullo Ismatov}

\address[N. J. Alves and H. Ismatov]{Applied Mathematics and Computational Sciences (AMCS), Computer, Electrical and Mathematical Sciences and Engineering Division (CEMSE), King Abdullah University of Science and Technology (KAUST), Thuwal, 23955-6900, Kingdom of Saudi Arabia}
\email[]{nuno.januarioalves@kaust.edu.sa}
\email[]{hikmatullo.ismatov@kaust.edu.sa}

\begin{document}

\begin{abstract}
We characterize the pointwise notion of almost equicontinuity for families of real-valued measurable functions on subsets of $\mathbb R^n$ of finite measure. The characterization is given by means of an integral truncated translation
condition. We also provide examples showing that the finite measure assumption and the truncation are essential.
\end{abstract}

\subjclass[2020]{Primary 28A20; Secondary 46E30, 54D30}
\keywords{Almost equicontinuity, almost equiboundedness, measurable functions,
convergence in measure, compactness, translation estimates}
\maketitle
\thispagestyle{empty}

\section{Introduction and main result}

Let $n\in\N$ and let $E$ be a measurable subset of $\R^n$. We consider families $\F$ of real-valued measurable functions defined on $E$. Throughout the paper, each $f\in\F$ is identified with its zero extension to $\R^n$, and for $y\in\R^n$ we denote by $\tau_y f$ the translation $\tau_y f(x):=f(x+y)$. Given a measurable set $A\subseteq\R^n$, we write $|A|$ for its Lebesgue measure.

Such a family $\F$ is said to be \emph{almost equicontinuous on} $E$ if, for each $\varepsilon>0$, there exists $\delta>0$ such that, for every $f\in\F$, there is a measurable set $B_f\subseteq E$ with $|B_f|<\varepsilon$ and
\[
|f(x_1)-f(x_2)|<\varepsilon
\]
whenever $x_1,x_2\in E\setminus B_f$ satisfy $|x_1-x_2|<\delta$. When $E$ has finite measure, this notion, together with almost equiboundedness, characterizes the relatively compact families of real-valued measurable functions on $E$ equipped with the topology of convergence in measure; see
\cite{hanson1933note,krotov2012criteria}. Related compactness criteria for spaces of measurable functions can also be found in~\cite{brudnyi2015compactness}. Recall that $\F$ is said to be \emph{almost equibounded on}  $E$ if, for each $\varepsilon>0$, there exists $M>0$ such that, for every $f\in\F$, there is a measurable set $S_f\subseteq E$ with $|S_f|<\varepsilon$ and
\[
|f|\le M
\qquad\text{on } E\setminus S_f.
\]

Almost equiboundedness has also appeared recently in the characterization of relatively compact subsets of asymptotic $L_p$ spaces on $\R^n$, $1\le p<\infty$. These spaces are denoted by $\Lambda^p(\R^n)$ and are defined
by
\begin{equation*} 
\Lambda^p(\R^n)
=
\left\{f:\R^n\to\R\text{ measurable}\ \Big| \
\int_{\R^n}\min(|f|,1)^p\,\dx<\infty\right\}.
\end{equation*}
For our purposes, it is enough to consider the case $p=1$. We recall that~$\Lambda^p(\R^n)$ is a completely metrizable topological vector space. We refer to~\cite{alves2025F} for the original introduction of the spaces $\Lambda_p$ as a class of almost-$L_p$ functions endowed with the topology of asymptotic $L_p$-convergence, to~\cite{alves2026truncation} for further developments including the equivalence between the original definition and the one considered here, and to~\cite{alves2026L1} for an application to $p$-Schr\"odinger equations.

In \cite{alves2026kolmogorov}, it has been proved that a family $\F\subseteq \Lambda^1(\R^n)$ is relatively compact if and only if it is almost equibounded, satisfies the truncated translation condition
\begin{equation} \label{eq:trunc_trans_int}
\lim_{|y|\to0}\sup_{f\in\F}
\int_{\R^n}\min(|f(x+y)-f(x)|,1)\,\dx=0,
\end{equation}
and satisfies the truncated tail condition
\begin{equation} \label{eq:trunc_tail}
\lim_{R\to\infty}\sup_{f\in\F}
\int_{\{|x|>R\}}\min(|f(x)|,1)\,\dx=0.
\end{equation}
An alternative proof of this compactness result has been obtained in~\cite{alves2026truncation} as a consequence of a truncation compactness criterion for asymptotic $L_p$ spaces on general measure spaces. Conditions~\eqref{eq:trunc_trans_int} and~\eqref{eq:trunc_tail} are the truncated counterparts of the translation and tail conditions in the classical Kolmogorov--Riesz theorem for $L^p(\R^n)$; see \cite{hanche2010kolmogorov, hanche2019improvement}. For related compactness
criteria in Banach function spaces, see~\cite{gorka2016arzela}.

We now make two observations. First, if all functions in the family $\F$ are supported in a common measurable set $E$ of finite measure, then condition~\eqref{eq:trunc_tail} is automatically satisfied. Indeed,
\[
\int_{\{|x|>R\}}\min(|f(x)|,1)\,\dx
\le
\big|\{|x|>R\}\cap E\big|
\to0
\qquad\text{as } R\to\infty,
\]
uniformly in $f\in\F$. Thus, for families supported in such a set $E$, the compactness criterion in $\Lambda^1(\R^n)$ reduces to almost equiboundedness together with condition~\eqref{eq:trunc_trans_int}.

Second, when $|E|<\infty$, every real-valued measurable function on $E$, identified with its zero extension to $\R^n$, belongs to $\Lambda^1(\R^n)$, since
\[
\int_{\R^n}\min(|f|,1)\,\dx
=
\int_E\min(|f|,1)\,\dx
\le |E|<\infty.
\]
Moreover, for two such functions $f$ and $g$, the distance induced by
$\Lambda^1(\R^n)$ is
\[
(f,g) \mapsto \int_E\min(|f-g|,1)\,\dx.
\]
This distance metrizes convergence in measure on $E$; see, for instance,~\cite{bogachev2007measure,alves2025F}. Hence, when $E$ has finite measure, we have, by the compactness theorem for
measurable functions recalled at the beginning, that relative compactness in~$\Lambda^1(\R^n)$ for families supported on $E$ amounts to almost equiboundedness and almost equicontinuity.

Combining these observations, one is naturally led to the conclusion that, for almost equibounded families of measurable functions on a set $E\subseteq\R^n$ of finite measure, almost equicontinuity should be equivalent to~\eqref{eq:trunc_trans_int}. The main point of this note is that the almost equiboundedness assumption is not needed. Moreover, one implication holds on arbitrary measurable subsets of $\R^n$. More precisely, we have:

\begin{theorem} \label{thm:main}
Let $E\subseteq\R^n$ be measurable, and let $\F$ be a family of real-valued measurable functions on $E$. Extend each $f\in\F$ by zero outside $E$. If $\F$ satisfies~\eqref{eq:trunc_trans_int}, then $\F$ is almost equicontinuous on $E$. Conversely, if $|E|<\infty$ and $\F$ is almost equicontinuous on $E$, then~\eqref{eq:trunc_trans_int} holds.
\end{theorem}

The paper is organized as follows. The implication from almost equicontinuity to~\eqref{eq:trunc_trans_int} is proved in Section~\ref{sec:necessity}; this is where the assumption $|E|<\infty$ is used. The implication from~\eqref{eq:trunc_trans_int} to almost equicontinuity is proved in Section~\ref{sec:sufficiency} and holds for every measurable set $E$. Finally, in Section~\ref{sec:examples} we present some examples. In particular, the first example shows that the finite measure assumption cannot be removed from the converse part of Theorem~\ref{thm:main}.

\section{Proof of Theorem~\ref{thm:main}: Necessity} \label{sec:necessity}

Let $\mathcal F$ be a family of real-valued measurable functions defined on a
measurable set $E\subseteq\R^n$ of finite measure. We extend each $f\in\mathcal F$ by zero outside $E$, and assume that $\mathcal F$ is almost equicontinuous on $E$.

Let $\eps>0$, and choose $\teps>0$ such that
\[
(3+|E|)\teps<\eps.
\]
By almost equicontinuity, there exists $\delta>0$ such that, for every
$f\in\mathcal F$, there is a measurable set $B_f\subseteq E$ with
$|B_f|<\teps$ and
\[
|f(x_1)-f(x_2)|<\teps
\]
whenever $x_1,x_2\in E\setminus B_f$ and $|x_1-x_2|<\delta$.

Since the functions are extended by zero outside $E$, we have, for every
$f\in\mathcal F$ and every $y\in\R^n$,
\[
|f(x+y)-f(x)|=0
\qquad
\text{for } x\in \bigl(E\cup(E-y)\bigr)^c,
\]
where
\[
E-y=\big\{x\in\R^n:\, x+y\in E\big\}.
\]

We write
\[
E\cup(E-y)=\bigl(E\cap(E-y)\bigr)\cup\bigl(E\triangle(E-y)\bigr),
\]
where $\triangle$ denotes the symmetric difference. Moreover,
\[
|E\triangle(E-y)|
=
\|\chi_E-\tau_y\chi_E\|_{L^1(\R^n)},
\]
where, for a measurable function $g$, we write $(\tau_y g)(x)=g(x+y)$. Since
$|E|<\infty$, we have $\chi_E\in L^1(\R^n)$. Hence, by the standard
continuity of translations in $L^1(\R^n)$; see, e.g.,
\cite[Proposition~8.5]{Folland}, it follows that
\[
|E\triangle(E-y)|\to0
\qquad\text{as } |y|\to0.
\]
Choose $r_0>0$ such that
\[
|E\triangle(E-y)|<\teps
\qquad\text{whenever } |y|<r_0,
\]
and set
\[
r:=\min\{\delta,r_0\}.
\]

Fix $f\in\mathcal F$ and $y\in\R^n$ with $|y|<r$. Then
\begin{align*}
\int_{\R^n}\min(|\tau_y f-f|,1)\,\dx
&=
\int_{E\cup(E-y)}\min(|\tau_y f-f|,1)\,\dx \\
&=
\int_{E\cap(E-y)}\min(|\tau_y f-f|,1)\,\dx \\
&\quad+
\int_{E\triangle(E-y)}\min(|\tau_y f-f|,1)\,\dx \\
&\le
\int_{E\cap(E-y)}\min(|\tau_y f-f|,1)\,\dx
+
|E\triangle(E-y)| \\
&<
\int_{E\cap(E-y)}\min(|\tau_y f-f|,1)\,\dx+\teps.
\end{align*}

It remains to estimate the integral over $E\cap(E-y)$. Set
\[
\widetilde B:=B_f\cup(B_f-y).
\]
Then $\widetilde B\subseteq E\cup(E-y)$ and
\[
|\widetilde B|\le 2|B_f|<2\teps.
\]
Therefore,
\begin{align*}
\int_{E\cap(E-y)}\min(|\tau_y f-f|,1)\,\dx
&=
\int_{(E\cap(E-y))\cap\widetilde B}
\min(|\tau_y f-f|,1)\,\dx \\
&\quad+
\int_{(E\cap(E-y))\setminus \widetilde B}
\min(|\tau_y f-f|,1)\,\dx \\
&\le
2\teps
+
\int_{(E\setminus B_f)\cap((E-y)\setminus(B_f-y))}
\min(|\tau_y f-f|,1)\,\dx.
\end{align*}

If
\[
x\in (E\setminus B_f)\cap((E-y)\setminus(B_f-y)),
\]
then $x,x+y\in E\setminus B_f$. Since
\[
|x+y-x|=|y|<r\le\delta,
\]
almost equicontinuity gives
\[
|f(x+y)-f(x)|<\teps.
\]
Hence
\[
\min(|f(x+y)-f(x)|,1)\le |f(x+y)-f(x)|<\teps
\]
on this set. Since this set is contained in $E$, we obtain
\[
\int_{(E\setminus B_f)\cap((E-y)\setminus(B_f-y))}
\min(|\tau_y f-f|,1)\,\dx
\le |E|\teps.
\]
Combining the estimates above yields
\[
\int_{\R^n}\min(|\tau_y f-f|,1)\,\dx
<
(3+|E|)\teps
<
\eps.
\]
This proves~\eqref{eq:trunc_trans_int}.
\qed

\section{Proof of Theorem~\ref{thm:main}: Sufficiency}
\label{sec:sufficiency}

Let $E\subseteq \R^n$ be a measurable set, and let $\mathcal F$ be a family of
real-valued measurable functions on $E$. We extend each $f\in\mathcal F$ by
zero outside $E$, and assume that \eqref{eq:trunc_trans_int} holds.

Fix $\varepsilon>0$, and set
\[
\lambda:=\min\left\{\frac{\varepsilon}{3},\frac12\right\}.
\]
By \eqref{eq:trunc_trans_int}, there exists $\delta>0$ such that
\[
\int_{\R^n}\min\big(|f(x+y)-f(x)|,1\big)\,\dx
<
\frac{\lambda\varepsilon}{2^{n+1}}
\]
for every $y\in\R^n$ with $|y|<\delta$ and every $f\in\mathcal F$.

Fix $f\in\mathcal F$. For $x\in E$, define
\[
A_f(x):=
\big\{z\in B(x,\delta):\, |f(z)-f(x)|\ge \lambda\big\},
\]
where $B(x,\delta)$ denotes the open ball centered at $x$ with radius $\delta$. Using the change of variables $z = x + y$ we obtain
\begin{align*}
|A_f(x)|  = \int_{\R^n}\chi_{A_f(x)}(z)\,\dz = \int_{\R^n}\chi_{A_f(x)}(x+y)\,\rd y.
\end{align*}
Now,
\[
\chi_{A_f(x)}(x+y)=1
\]
if and only if
\[
x+y\in B(x,\delta)
\quad\text{and}\quad
|f(x+y)-f(x)|\ge \lambda.
\]
Since $x+y\in B(x,\delta)$ is equivalent to $y\in B(0,\delta)$, it follows that
\[
|A_f(x)|
=
\int_{B(0,\delta)}
\chi_{\{|f(x+y)-f(x)|\ge \lambda\}}\,\rd y.
\]

Integrating over $E$ and using Tonelli's theorem, we get
\[
\begin{aligned}
\int_E |A_f(x)|\,\dx
&=
\int_E\int_{B(0,\delta)}
\chi_{\{|f(x+y)-f(x)|\ge \lambda\}}\,\rd y\,\dx  \\
&=
\int_{B(0,\delta)}
\int_E
\chi_{\{|f(x+y)-f(x)|\ge \lambda\}}\,\dx\,\rd y  \\
&=
\int_{B(0,\delta)}
\big|\{x\in E:\ |f(x+y)-f(x)|\ge \lambda\}\big|\,\rd y.
\end{aligned}
\]
Since $E\subseteq\R^n$, we have
\[
\big|\{x\in E:\ |f(x+y)-f(x)|\ge \lambda\}\big|
\le
\big|\{x\in \R^n:\ |f(x+y)-f(x)|\ge \lambda\}\big|.
\]
Moreover, because $0<\lambda\le 1/2$, for every real number $u$ one has
\[
\lambda\chi_{\{|u|\ge \lambda\}}\le \min(|u|,1).
\]
Applying this with $u=f(x+y)-f(x)$ gives, for every $|y|<\delta$,
\[
\begin{aligned}
\big|\{x\in \R^n:\ |f(x+y)-f(x)|\ge \lambda\}\big|
&\le
\frac1\lambda
\int_{\R^n}\min\big(|f(x+y)-f(x)|,1\big)\,\dx  \\
&<
\frac{\varepsilon}{2^{n+1}}.
\end{aligned}
\]
Consequently,
\[
\int_E |A_f(x)|\,\dx
<
\frac{\varepsilon}{2^{n+1}}\,|B(0,\delta)|.
\]

Now define
\[
B_f:=
\left\{x\in E:\ |A_f(x)|>\frac{|B(0,\delta)|}{2^{n+1}}\right\},
\]
and note that
\[
|B_f|
\le
\frac{2^{n+1}}{|B(0,\delta)|}
\int_E |A_f(x)|\,\dx
<
\varepsilon.
\]

We show that $\delta$ gives the desired modulus of almost
equicontinuity. Let $x_1,x_2\in E\setminus B_f$ be such that
\[
|x_1-x_2|<\delta.
\]
Since $x_j\notin B_f$, $j=1,2$, we have
\[
|A_f(x_j)|\le \frac{|B(0,\delta)|}{2^{n+1}},
\qquad j=1,2.
\]
Set
\[
I:=B(x_1,\delta)\cap B(x_2,\delta).
\]
The set $I$ contains the ball
\[
B\left(\frac{x_1+x_2}{2},\,\delta-\frac{|x_1-x_2|}{2}\right).
\]
Since $|x_1-x_2|<\delta$, this ball has radius strictly larger than
$\delta/2$. Hence
\[
|I|>2^{-n}|B(0,\delta)|.
\]
On the other hand,
\[
\begin{aligned}
|A_f(x_1)\cup A_f(x_2)|
&\le
|A_f(x_1)|+|A_f(x_2)|  \\
&\le
2^{-n}|B(0,\delta)|.
\end{aligned}
\]
Therefore
\[
|I|>|A_f(x_1)\cup A_f(x_2)|.
\]
Now we claim that this implies
\[
I\setminus\big(A_f(x_1)\cup A_f(x_2)\big)\neq\varnothing.
\]
Indeed, if
\[
I\setminus\big(A_f(x_1)\cup A_f(x_2)\big)=\varnothing,
\]
then
\[
I\subseteq A_f(x_1)\cup A_f(x_2),
\]
and hence
\[
|I|\le |A_f(x_1)\cup A_f(x_2)|,
\]
which contradicts the strict inequality above. Thus the set difference is
nonempty.

Choose
\[
z\in I\setminus\big(A_f(x_1)\cup A_f(x_2)\big).
\]
Since $z\in I$, we have
\[
z\in B(x_1,\delta)
\quad\text{and}\quad
z\in B(x_2,\delta).
\]
Since also
\[
z\notin A_f(x_1)\cup A_f(x_2),
\]
we have $z\notin A_f(x_1)$ and $z\notin A_f(x_2)$. By the definition of
$A_f(x)$, this gives
\[
|f(z)-f(x_1)|<\lambda
\quad\text{and}\quad
|f(z)-f(x_2)|<\lambda.
\]
Consequently,
\[
|f(x_1)-f(x_2)|
\le
|f(x_1)-f(z)|+|f(z)-f(x_2)|
<
2\lambda
<
\varepsilon.
\]

Thus, for the given $\varepsilon>0$, we have found $\delta>0$ such that for
every $f\in\mathcal F$ there exists a measurable set $B_f\subseteq E$ with
$|B_f|<\varepsilon$ and
\[
|f(x_1)-f(x_2)|<\varepsilon
\]
whenever $x_1,x_2\in E\setminus B_f$ and $|x_1-x_2|<\delta$. This proves that
$\mathcal F$ is almost equicontinuous.
\qed

\section{Examples} \label{sec:examples}

We begin with an example showing that the assumption $|E|<\infty$ is indispensable for the implication from almost equicontinuity to~\eqref{eq:trunc_trans_int}.

\begin{example}\label{ex:finite-measure-needed}
Let $E=\R$, and for $m\in\N$ define
\[
f_m(x):=x \, \chi_{[0,m]}(x),\qquad x\in\R.
\]
Then the family
\[
\mathcal F:=\big\{f_m:\, m\in\N\big\}
\]
is almost equicontinuous on $\R$, but it does not satisfy~\eqref{eq:trunc_trans_int}.

Indeed, fix $\varepsilon>0$, and choose
\[
0<\delta<\frac{\varepsilon}{4}.
\]
For each $m\in\N$, set
\[
B_m:=(m-\delta,m+\delta), \qquad |B_m|=2\delta<\varepsilon.
\]
If
$x_1,x_2\in\R\setminus B_m$ and $|x_1-x_2|<\delta$, then $x_1$ and $x_2$
cannot lie on opposite sides of the interval $B_m$. On $(-\infty,m-\delta]$,
the function $f_m$ agrees with $x\mapsto \max\{x,0\}$, which is Lipschitz
with constant $1$, while on $[m+\delta,\infty)$ the function $f_m$ vanishes.
Hence
\[
|f_m(x_1)-f_m(x_2)|<\varepsilon.
\]
Thus $\mathcal F$ is almost equicontinuous on $\R$.

On the other hand, let $0<y<1$. Since $m>y$ for every $m\in\N$, for
each $m\in\N$ and every $x\in[0,m-y]$ we have $x,x+y\in[0,m]$, and therefore
\[
f_m(x+y)-f_m(x)=y.
\]
Consequently,
\[
\int_{\R}\min(|f_m(x+y)-f_m(x)|,1)\,\dx
\ge
\int_0^{m-y} y\,\dx
=
y(m-y).
\]
Letting $m\to\infty$, we obtain
\[
\sup_{m\in\N}
\int_{\R}\min(|f_m(x+y)-f_m(x)|,1)\,\dx
=
\infty
\]
for every $0<y<1$. Therefore condition~\eqref{eq:trunc_trans_int} does not hold.
\end{example}

The next example shows that the truncation in \eqref{eq:trunc_trans_int} is
essential. Even on a set of finite measure, almost equicontinuity need not
imply the untruncated translation condition
\begin{equation} \label{eq:untrunc}
\lim_{|y|\to0}\sup_{f\in\mathcal F}
\int_{\R^n}|f(x+y)-f(x)|\,\dx=0.
\end{equation}

\begin{example}\label{ex:truncation-is-essential}
Let $E=(0,1)$ and, for $m\in\N$, set
\[
I_m:=\left(0,\frac1m\right),
\qquad
f_m:=m \, \chi_{I_m}.
\]
Let
\[
\mathcal F:=\big\{f_m:m\in\N\big\}.
\]
We extend each $f_m$ by zero outside $E$.

First, $\mathcal F$ satisfies the truncated translation condition. Indeed,
since
\[
\tau_y f_m=m \, \chi_{I_m-y},
\]
we have
\[
\min(|\tau_y f_m-f_m|,1)
=
\chi_{I_m\triangle(I_m-y)}.
\]
Therefore
\[
\int_{\R}\min(|\tau_y f_m-f_m|,1)\,\dx
=
|I_m\triangle(I_m-y)|
\le
2|y|.
\]
Consequently,
\[
\sup_{m\in\N}
\int_{\R}\min(|\tau_y f_m-f_m|,1)\,\dx
\le
2|y|\to0
\qquad\text{as } y\to0.
\]
Thus \eqref{eq:trunc_trans_int} holds, and so $\mathcal F$ is almost
equicontinuous on $E$ by Theorem~\ref{thm:main}.

However, the untruncated translation condition fails. Fix $y\ne0$. If
$m>|y|^{-1}$, then $|y|>1/m$, and the intervals $I_m$ and $I_m-y$ are
disjoint. Hence
\[
\int_{\R}|\tau_y f_m-f_m|\,\dx
=
m|I_m\triangle(I_m-y)|
=
m\left(\frac2m\right)
=
2.
\]
Thus
\[
\sup_{m\in\N}
\int_{\R}|\tau_y f_m-f_m|\,\dx
\ge 2
\]
for every $y\ne0$. Therefore condition~\eqref{eq:untrunc} does not
hold.
\end{example}

In view of Theorem~\ref{thm:main}, it is natural to ask whether an analogous
statement holds for equicontinuity if one replaces the truncated condition
\eqref{eq:trunc_trans_int} by the untruncated condition~\eqref{eq:untrunc}.
The last example shows that this is not the case, even for equibounded families of continuous functions. 

Recall that a family $\mathcal G$ of real-valued functions on $E$ is
\emph{equicontinuous} if, for every $\varepsilon>0$, there exists $\delta>0$ such that
\[
|g(x_1)-g(x_2)|<\varepsilon
\]
for every $g\in\mathcal G$ and all $x_1,x_2\in E$ satisfying
$|x_1-x_2|<\delta$, and it is \emph{equibounded} if there exists $M>0$ such that
\[
|g(x)|\le M
\]
for every $g\in\mathcal G$ and every $x\in E$.

\begin{example}\label{ex:no-equicontinuity-criterion}
Let $E=(0,1)$. For $m\ge 3$, define $g_m:E\to\R$ by
\[
g_m(x)=
\begin{dcases}
0, & \text{if }0<x\le \dfrac12, \\
m\left(x-\dfrac12\right), & \text{if } \dfrac12<x<\dfrac12+\dfrac1m, \\
1, & \text{if } \frac12+\frac1m\le x<1.
\end{dcases}
\]
Then each $g_m$ is continuous on $(0,1)$ and
\[
0\le g_m\le 1.
\]
Thus the family
\[
\mathcal G:=\big\{g_m:\, m\ge 3\big\}
\]
is equibounded.

However, $\mathcal G$ is not equicontinuous. Indeed, given $\delta>0$, choose
$m>\max\{3,\delta^{-1}\}$, and set
\[
x_1=\frac12,
\qquad
x_2=\frac12+\frac1m.
\]
Then $x_1,x_2\in(0,1)$,
\[
|x_1-x_2|=\frac1m<\delta,
\]
but
\[
|g_m(x_1)-g_m(x_2)|=1.
\]

On the other hand, the family satisfies condition~\eqref{eq:untrunc}. Extend each $g_m$ by zero outside~$(0,1)$. We claim that
\[
\int_{\R}|\tau_y g_m-g_m|\,\dx\le 2|y|
\]
for every $m\ge3$ and every $y\in\R$ with $|y|<1/4$. 

It is enough to consider $0<y<1/4$, since the case $y<0$ follows by a change
of variables. The zero extension of $g_m$ vanishes on $(-\infty,1/2]$, and
$g_m$ is nondecreasing on $(0,1)$. Using these facts and the assumption
$0<y<1/4$, we obtain
\begin{align*}
\int_{\R}|g_m(x+y)-g_m(x)|\,\dx
&=
\int_0^{1-y}\bigl(g_m(x+y)-g_m(x)\bigr)\,\dx
+
\int_{1-y}^1 g_m(x)\,\dx  \\
&=
\int_{1-y}^1 g_m(x)\,\dx
-
\int_0^y g_m(x)\,\dx
+
\int_{1-y}^1 g_m(x)\,\dx.
\end{align*}
Moreover, because $0<y<1/4$, we have $g_m=0$ on $(0,y)$. Since also
$0\le g_m\le1$, it follows that
\[
\int_{\R}|g_m(x+y)-g_m(x)|\,\dx
\le
2\int_{1-y}^1 g_m(x)\,\dx
\le
2y.
\]
Therefore
\[
\sup_{m\ge3}
\int_{\R}|\tau_y g_m-g_m|\,\dx
\le
2|y|\to0
\qquad\text{as } y\to0.
\]
Thus $\mathcal G$ satisfies~\eqref{eq:untrunc}, but it is not equicontinuous.
\end{example}

\section*{Acknowledgments}

This publication is based upon work supported by King Abdullah University of Science and Technology (KAUST) under Award No. ORFS-CRG12-2024-6430.

\end{document}